%&amstex
\input amstex
\input amsppt.sty
\magnification=\magstep1
\hsize=30truecc
\vsize=22.2truecm
\baselineskip=16truept
\NoBlackBoxes
\TagsOnRight \pageno=1 \nologo
\def\Z{\Bbb Z}
\def\N{\Bbb N}

\def\l{\left}
\def\r{\right}
\def\bg{\bigg}
\def\({\bg(}
\def\[{\bg\lfloor}
\def\){\bg)}
\def\]{\bg\rfloor}
\def\t{\text}
\def\f{\frac}

\def\bi{\binom}
\def\eq{\equiv}

\def\ls{\leqslant}

\def\da{\delta}

\def\Proof{\noindent{\it Proof}}

\hbox {Preprint}
\bigskip
\topmatter
\title On convolutions of Euler numbers\endtitle
\author Zhi-Wei Sun\endauthor
\leftheadtext{Zhi-Wei Sun} \rightheadtext{On convolutions of
Euler numbers}
\affil Department of Mathematics, Nanjing University\\
 Nanjing 210093, People's Republic of China
  \\  zwsun\@nju.edu.cn
  \\ {\tt http://math.nju.edu.cn/$\sim$zwsun}
\endaffil
\abstract  We show that if $p$ is an odd prime then
$$\sum_{k=0}^{p-1}E_kE_{p-1-k}\eq1\pmod p$$
and
$$\sum_{k=0}^{p-3}E_kE_{p-3-k}\eq(-1)^{(p-1)/2}2E_{p-3}\pmod p,$$
where $E_0,E_1,E_2,\ldots$ are Euler numbers.
Moreover, we prove that for any positive integer $n$ and prime number $p>2n+1$ we have
$$\sum_{k=0}^{p-1+2n}E_kE_{p-1+2n-k}\eq s(n)\pmod p$$
where $s(n)$ is an integer only depending on $n$.
\endabstract
\thanks 2010 {\it Mathematics Subject Classification}.\,Primary 11B68;
Secondary 11A07.
\newline\indent Supported by the National Natural Science
Foundation (grant 10871087) and the Overseas Cooperation Fund (grant 10928101) of China.
\endthanks
\endtopmatter
\document

\heading{1. Introduction}\endheading

The Euler numbers $E_n\ (n\in\N=\{0,1,2,\ldots\})$ are integers
defined by
$$E_0=1\ \t{and}\ \sum^n\Sb k=0\\2\mid k\endSb \bi nk E_{n-k}=0\ \ \ \t{for}\ n\in\Z^+=\{1,2,3,\ldots\}.$$
It is well known that $E_{2n+1}=0$ for all $n\in\N$ and
$$\sec x=\sum_{n=0}^\infty(-1)^n E_{2n}\f{x^{2n}}{(2n)!}\ \ \l(|x|<\f{\pi}2\r).$$
The exponential generating function for Euler numbers is given by
$$\f{2e^x}{e^{2x}+1}=\sum_{n=0}^\infty E_n\f{x^n}{n!}\quad \ \l(|x|<\f{\pi}2\r).$$
Thus
$$\l(\f{2e^x}{e^{2x}+1}\r)^2=\sum_{k=0}^\infty E_k\f{x^k}{k!}\sum_{l=0}^\infty E_l\f{x^l}{l!}
=\sum_{n=0}^\infty f(n)\f{x^n}{n!},$$
where
$$f(n)=\sum_{k=0}^n\bi nk E_kE_{n-k}.$$
In this paper we are interested in the usual convolution of Euler numbers given by
$\sum_{k=0}^nE_kE_{n-k}$. The reader may consult [PS], [SP] and [S11] for related background.

Now we present our main results.

\proclaim{Theorem 1.1} Let $p$ be an odd prime. Then
$$\sum_{k=0}^{p-3}E_kE_{p-3-k}\eq2\l(\f{-1}p\r)E_{p-3}\pmod p,\tag1.1$$
where $(-)$ denotes the Jacobi symbol.
Moreover, for any $n=0,1,2,\ldots$ we have
$$\sum_{k=0}^{p-1+2n}E_kE_{p-1+2n-k}\eq s(n)+\da(p,n)\pmod p,\tag1.2$$
where
$$s(n)=\sum_{k=0}^nE_{2k}E_{2n-2k}\tag1.3$$
and
$$\da(p,n)=\cases1&\t{if}\ n>0\ \&\ p-1\mid 2n,
\\0&\t{otherwise}.\endcases$$
\endproclaim

{\it Example} 1.1. Here are the values of $s(n)$ with $n\in\{0,1,2,3,4,5\}$:
$$s(0)=1,\ s(1)=-2,\ s(2)=11,\ s(3)=-132,\ s(4)=2917,\ s(5)=-104422.$$
Thus, for any odd prime $p$ we have
$$\align \sum_{k=0}^{p-1}E_kE_{p-1-k}\eq&1\pmod p,\tag1.4
\\\sum_{k=0}^{p+1}E_kE_{p+1-k}\eq&-2\pmod p\ \ \t{if}\ p>3,\tag1.5
\\\sum_{k=0}^{p+3}E_kE_{p+3-k}\eq&11\pmod p\ \ \t{if}\ p>5.\tag1.6
\endalign$$

Applying (1.2) again and again we immediately obtain the following consequence.

\proclaim{Corollary 1.1} Let $n=\f{p-1}2q+r$ with $q\in\{1,2,3,\ldots\}$ and $r\in\{0,\ldots (p-3)/2\}$. Then we have
$$s(n)\eq s(r)+(q-1)\da_{r,0}\pmod p.\tag1.7$$
\endproclaim

By a further refinement of our method to prove Theorem 1.1 and some complicated discussions,
we can deduce the following theorem
though we will not give the details of the proof since it is similar to that of Theorem 1.1.

\proclaim{Theorem 1.2} For any odd prime $p$, we have
$$\sum_{i+j+k=p-3}E_iE_jE_k\eq-2E_{p-3}\pmod p.\tag1.8$$
Also, for each $n\in\N$ there is a unique integer $t(n)$ such that if $p>2n+1$ is a prime then
$$\sum_{i+j+k=p-1+2n}E_iE_jE_k\eq t(n)\pmod p.\tag1.9$$
In particular,
$$t(0)=3,\ t(1)=-9,\ t(2)=68,\ t(3)=-1068.$$
\endproclaim

Theorems 1.1 and 1.2 should have their $q$-analogues. We leave this to those who are interested in such things.

\heading{2. Proof of Theorem 1.1}\endheading

\proclaim{Lemma 2.1} Let $p$ be an odd prime and let $k\in\N$ be even. Then
$$E_k\eq2\sum^{p-1}\Sb j=1\\2\nmid j\endSb\l(\f{-1}j\r)j^k+\da_{k,0}\l(\f{-1}p\r)\pmod{p}.\tag2.1$$
\endproclaim
\Proof. By [S05, (1.1)],
$$E_k\eq\sum_{i=0}^{p-1}(-1)^i(2i+1)^k\pmod p.$$
Observe that
$$\align&\sum_{i=0}^{p-1}(-1)^i(2i+1)^k-(-1)^{(p-1)/2}p^k
\\=&\sum_{i=0}^{(p-3)/2}\l((-1)^i(2i+1)^k+(-1)^{p-1-i}(2(p-1-i)+1)^k\r)
\\\eq&2\sum_{i=0}^{(p-3)/2}(-1)^i(2i+1)^k=2\sum^{p-1}\Sb j=1\\2\nmid j\endSb\l(\f{-1}j\r)j^k\pmod p.
\endalign$$
So (2.1) follows. \qed

\medskip
\noindent{\it Proof of Theorem 1.1}.
(i) In view of Lemma 2.1,
$$\align\sum_{k=0}^{p-3}E_kE_{p-3-k}\eq&2\l(\f{-1}p\r)\times2\sum^{p-1}\Sb j=1\\2\nmid j\endSb\l(\f{-1}j\r)j^{p-3}
\\&+2\sum_{k=0}^{(p-3)/2}\sum^{p-1}\Sb i=1\\2\nmid i\endSb\l(\f{-1}i\r)i^{2k}2\sum^{p-1}\Sb j=1\\2\nmid j\endSb\l(\f{-1}j\r)j^{p-3-2k}
\\\eq&2\l(\f{-1}p\r)E_{p-3}+4\sum^{p-1}\Sb j=1\\2\nmid j\endSb\l(\f{-1}j\r)^2j^{p-3}
\\&+8\sum\Sb 1\ls i<j<p\\2\nmid ij\endSb\l(\f{-1}{ij}\r)j^{p-3}\f{(i^2/j^2)^{(p-1)/2}-1}{i^2/j^2-1}
\\\eq&2\l(\f{-1}p\r)E_{p-3}+4\sum_{j=1}^{p-1}\f1{j^2}-4\sum_{k=1}^{(p-1)/2}\f1{(2k)^2}
\\\eq&2\l(\f{-1}p\r)E_{p-3}\pmod p.
\endalign$$
In the last step we noted that
$$2\sum_{k=1}^{(p-1)/2}\f1{k^2}\eq\sum_{k=1}^{(p-1)/2}\l(\f1{k^2}+\f1{(p-k)^2}\r)=\sum_{k=1}^{p-1}\f1{k^2}\eq0\pmod p$$
by the Wolstenhomle congruence. Thus (1.1) holds.

(ii) Observe that
$$\sum_{k=0}^{p-1+2n}E_kE_{p-1+2n-k}=\sum_{k=0}^{(p-3)/2}E_{2k}E_{p-1+2n-2k}+\sum_{k=0}^nE_{p-1+2k}E_{2n-2k}.\tag2.2$$
By Lemma 2.1,
$$\align E_{p-1}\eq&2\sum^{p-1}\Sb j=1\\2\nmid j\endSb\l(\f{-1}j\r)
=2\l(1-1+\cdots+(-1)^{(p-3)/2}(p-2)\r)
\\=&1-\l(\f{-1}p\r)=E_0-\l(\f{-1}p\r)\pmod{p}
\endalign$$
and also
$$E_{p-1+2k}\eq E_{2k}\pmod p\quad\t{for}\ k=1,2,3,\ldots.$$
Therefore
$$\sum_{k=0}^nE_{p-1+2k}E_{2n-2k}\eq\sum_{k=0}^nE_{2k}E_{2n-2k}-\l(\f{-1}p\r)E_{2n}\pmod p.\tag2.3$$

In view of Lemma 2.1, we also have
$$\align &\sum_{k=0}^{(p-3)/2}E_{2k}E_{p-1+2n-2k}
\\\eq&\l(\f{-1}p\r)2\sum^{p-1}\Sb j=1\\2\nmid j\endSb\l(\f{-1}j\r)j^{p-1+2n}
\\&+\sum_{k=0}^{(p-3)/2}2\sum\Sb 0<i<p\\2\nmid i\endSb\l(\f{-1}i\r)i^{2k}2\sum\Sb 0<j<p\\2\nmid j\endSb\l(\f{-1}j\r)j^{p-1+2n-2k}
\\\eq&\l(\f{-1}p\r)\l(E_{2n}-\da_{n,0}\l(\f{-1}p\r)\r)+4\sum\Sb 0<i,j<p\\2\nmid ij\endSb\l(\f{-1}{ij}\r)j^{2n}\sum_{k=0}^{(p-3)/2}\f {i^{2k}}{j^{2k}}
\\\eq&\l(\f{-1}p\r)E_{2n}-\da_{n,0}+4\sum\Sb 0<j<p\\2\nmid j\endSb\l(\f{-1}{j^2}\r)j^{2n}\f{p-1}2
\\&+8\sum\Sb 0<i<j<p\\2\nmid ij\endSb\l(\f{-1}{ij}\r)j^{2n}\f{(i^2/j^2)^{(p-1)/2}-1}{i^2/j^2-1}
\endalign$$
Note that if $0<i,j<p$ and $2\nmid ij$ then $i\not\eq -j\pmod p$ since $i+j$ is even while $p$ is odd.
Applying Fermat's little theorem we obtain from the above
$$\sum_{k=0}^{(p-3)/2}E_{2k}E_{p-1+2n-2k}\eq
\l(\f{-1}p\r)E_{2n}-\da_{n,0}-2\sum\Sb 0<j<p\\2\nmid j\endSb j^{2n}\pmod{p}.$$
If $p-1$ divides $2n$, then
$$\sum\Sb 0<j<p\\2\nmid j\endSb j^{2n}\eq|\{0<j<p:\ 2\nmid j\}|=\f{p-1}2\pmod p.$$
When $p-1\nmid 2n$, we have
$$2\sum_{j=1}^{(p-1)/2}j^{2n}\eq\sum_{j=1}^{(p-1)/2}\l(j^{2n}+(p-j)^{2n}\r)=\sum_{j=1}^{p-1}j^{2n}\eq0\pmod p$$
(cf. [IR, p.\,235])
and hence
$$\sum\Sb 0<j<p\\2\nmid j\endSb j^{2n}=\sum_{j=1}^{p-1}j^{2n}-\sum_{j=1}^{(p-1)/2}(2j)^{2n}\eq0\pmod p.$$
Thus
$$\sum_{k=0}^{(p-3)/2}E_{2k}E_{p-1+2n-2k}\eq\l(\f{-1}p\r)E_{2n}-\da_{n,0}+[p-1\mid 2n]\pmod p,\tag2.4$$
where $[p-1\mid 2n]$ takes $1$ or $0$ according as $p-1\mid 2n$ or not.

Combining (2.2)-(2.4) we get
$$\sum_{k=0}^{p-1+2n}E_kE_{p-1+2n-k}\eq\sum_{k=0}^nE_{2k}E_{2n-2k}-\da_{n,0}+[p-1\mid 2n]\pmod{p}.$$
This proves (1.2).

So far we have completed the proof of Theorem 1.1. \qed

\enddocument

\bye